\documentclass{article}
\usepackage{amssymb}
\newcommand {\ex} {\mathbb{E}}
\newcommand {\vs} {\varsigma}
\newcommand {\vt} {\vartheta}
\newcommand {\cL} {\mathcal{L}}
\newcommand {\cX} {\mathcal{X}}

\newcommand {\syF} {\mbox{\sy F}}
\newcommand {\syG} {\mbox{\sy G}}
\newcommand {\syH} {\mbox{\sy H}}

\newcommand {\fr} {\mbox{\bf r}}
\newcommand {\fs} {\mbox{\bf s}}

\newcommand {\fone} {\mbox{\bf 1}}
\newcommand {\st} {\scriptscriptstyle}
\newtheorem{lemma}{Lemma}
\newtheorem{theorem}{Theorem}
\newtheorem{prop}{Proposition}
\newcommand{\reff}[1]{(\ref{#1})}
\newfont{\sy}{cmsy10 scaled 1000}
\begin{document}
\title{Proliferation Model Dependence in Fluctuation Analysis: The Neutral Case}
\date{}
\author{Wolfgang P. Angerer\\Goethe-Universit\"at\\Frankfurt am Main\\Germany}
\maketitle
\newpage
\begin{abstract}
\noindent We discuss the evaluation of Luria-Delbr\"uck fluctuation experiments under Bellman-Harris models of cell proliferation. It is shown that under certain very natural assumptions concerning the life-time distributions and the offspring distributions of mutant and non-mutant bacteria, the suitably normed and centered number of mutants contained in a large culture of bacteria converges to a stable random variable with index 1. The result obtains under the assumption that the mutation under consideration is \lq neutral' in the sense that on average, mutant cells produce the same number of offspring as non-mutant cells. Thus, provided this condition is met, fluctuation experiments can be evaluated without knowledge of how cell proliferation proceeds in detail. This resolves a long-standing issue in the theory of Luria-Delbr\"uck fluctuation analysis.\\


\end{abstract}
\newpage
\section{Introduction} \label{intro}
Fluctuation analysis, or the art of measuring mutation rates by means of Luria and Delbr\"uck's fluctuation test, has come into prominence sixty-odd years ago \cite{ld}, when it was realized that the huge fluctuations of the number of phage-resistant bacteria in sensitive cultures of {\it Escherichia coli} could be taken as evidence that the phage-resistant bacteria had arisen in these cultures for no particular reason: For if they had, what better reason could there be for a bacterium to become phage-resistant than the encounter with the phage itself? But then comparable cultures of bacteria should, upon exposure to comparable amounts of the phage, produce more or less the same number of mutant colonies on solid medium. They don't, even if one allows for some variability in the number of mutant colonies a given culture of bacteria produces in an actual experiment, and so the development of phage-resistance by these bacteria is not triggered by the phage. Indeed, under the assumption that resistant bacteria are only produced in the presence of the phage, the actual number of mutant colonies emerging from reasonably many cultures of bacteria on solid medium should be distributed over the cultures in a Poisson-like manner. This is not the case, and so there must be at work some other mechanism by which sensitive bacteria are turned into resistant ones. Fluctuation analysis aims to illuminate this mechanism by way of a thorough analysis of the distribution of mutant colony count in what has come to be called a {\it fluctuation experiment}, which is a fancy term for the fact that a number of bacterial cultures have been grown from an initially  small number of non-mutant bacteria, and the number of mutants each culture contains has been recorded (see, e.g., \cite{wa1, as, btn, kf, ko, ka, ok, tb} for more or less recent reviews of how fluctuation experiments are performed in the laboratory). A particular emphasis of this analysis lies on determining the probability $\rho$ (not a standard notation) for a cell to undergo mutation during division. It is customary to refer to this probability as the {\it mutation rate}.

It is true that the fluctuation-analytic approach has faced recent competition from certain other, experimentally more sophisticated methods \cite{rf, bd, f, ngr1, ngr2, rgn}. Still, it may serve as a welcome complement to these, all the more as certain parameters of possible interest besides the mutation rate can be determined from a fluctuation experiment as well \cite{wa2}. For this, one requires a mathematical model of how to calculate the so-called {\it Luria-Delbr\"uck distribution} of mutant cells an individual culture of bacteria contains just before it is plated on solid medium. Now, it is not only since the (relatively recent) induction of fluctuation analysis into the \lq molecular biology wing of the Museum of Elegant Science' (the copyright for this {\it nice} one goes to Patricia Foster \cite{f}) that questions concerning the mathematical theory of bacterial populations becoming enriched with mutants have been paid a certain amount of attention. Not a very big lot of it, mind you, but the occasional article on the theory of the Luria-Delbr\"uck distribution has popped up in the journals for longer than the recent half of a century (e.g., \cite{kemp, kend, lc, mss, pakes, p, sar, sms, sgl, qz1, qz2, qz3} besides those already quoted), and will, in all likelihood, continue to do so until the community is agreed that the issue simply isn't interesting anymore. To the author's knowledge, however, the influence that different models of cell proliferation might exert on the calculation of the Luria-Delbr\"uck distribution and the evaluation of a fluctuation experiment has never been touched upon in detail. What one needs to look at in particular is {\it non-exponential} life-time distributions of cells. This is what makes clonal expansion (if only cell numbers are considered and not, say, the distribution of cells over possibly a continuum of age-classes) non-Markov and, therefore, a little messy. True enough, the probability generating function (PGF) of the Luria-Delbr\"uck distribution is, the non-Markovian character of bacterial proliferation notwithstanding, easy to write down. But what's next? Even with the help of a computer program for symbolic calculation, the actual computation of the PGF in question isn't all that easy. (The reader who disagrees is cordially invited to try her hand at Equation \reff{basic} below  - and no, it is {\it not} permissible to assume that life-time distributions are exponential!) Is it all worth it?

It could be argued that, to be on the safe side, one could always conduct a fluctuation experiment such that a reasonable fraction $p_0$ of cultures will not contain any mutants at all. Since the probability that an initially very small culture of bacteria will not have produced {\it any} mutants by the time it has reached a certain size $n$ (which is to say, after $\sim n$ cell divisions) is $(1 - \rho)^n$, an obvious estimate for the mutation rate would be $\rho \sim - \log p_0/n$. This would be true regardless of how the proliferation of cells proceeds in detail. The performance of this so-called $p_0$-method, however, is generally poor \cite{lc}, so that the question is whether other methods of evaluating a fluctuation experiment need to be adapted to non-Markovian  models of cell proliferation. We will see that even for \lq neutral' models of mutation, the answer is in general yes. Just how much these methods need to be adapted is the main result of this paper. A translation of what should be understood by \lq neutrality' in the context of single cell populations is as follows: a mutation is \lq neutral' if, on average and in the long run, a mutant cell produces as many offspring as a non-mutant one. We shall be more precise about this in the following section.
\section{The Luria-Delbr\"uck Distribution}
We shall take it for granted that the PGF of the Luria-Delbr\"uck distribution is of the form
\begin{equation}
g_{\st LD}(s) :=: \sum_{r=0}^{\infty} p_r s^r =  \exp\Big(m\big(g(s)-1\big)\Big)\:,\label{basic}
\end{equation}
where $m$ is parameter that we shall soon return to,
\begin{equation}
g(s) := \beta \int_0^{\infty} e^{- \beta u} F_u(s) \, du \:, \label{g}
\end{equation}
and
\begin{equation}
F_u (s) = \big(1 - G^{\bullet}(u)\big)s + \int_0^u f \circ F_{u - y}(s) \, dG^{\bullet}(y) \label{iq}
\end{equation}
is the PGF of particle numbers in a Bellman-Harris process with first-generation offspring PGF
\begin{equation}
f(s) := \sum_{k=0}^{\infty} \pi_k s^k\:,\label{f}
\end{equation}
and life-time distribution function $G^{\bullet}(\cdot)$ \cite{an}. In words, $G^{\bullet}(t)$ is the probability that a newborn mutant cell lives at least until time $t$, and $\pi_k$ is the probability that once it splits into progeny, it will split into exactly $k$ of these. We will always assume that proliferation is supercritical, that is, $1 < f'(1) =: \mu^{\bullet} < \infty$. This implies, in particular, that there exists $q \in [0,1)$ such that $f(q) = q$. For practical purposes, of course, $f(s)$ can be taken to be only a quadratic polynomial in $s$, but this doesn't make anything about the calculations easier. We shall give little justification for Equation \reff{basic} in this paper; it is in agreement with (and in fact a generalization of) previous work on the subject (notably \cite{lc, sgl, tan1, tan2, dlm}), and in any case seems reasonable enough: the point is that with the mutation rate $\rho$ being generally small, the first mutants will arise only when the population as a whole is already well into the stage of Malthusian growth, so that the probability for a mutant to arise $u$ units of time before the present moment is $\sim e^{-\beta u}$ smaller. Here and in Equation \reff{g}, $\beta$ is the Malthusian parameter for the growth of a population of {\it exclusively} non-mutant cells, which for our purposes is a number that satisfies
\begin{equation}
\mu^{\circ} \, \int_0^{\infty} e^{- \beta t} \, dG^{\circ}(t) = 1 \:,\label{nmmalthus}
\end{equation}
where $\mu^{\circ}$ is the expected number of non-mutant first-generation offspring of a non-mutant cell, and $G^{\circ}(\cdot)$ is the life-time distribution of non-mutant bacteria. We now define a mutation to be \lq neutral' if $\beta$ also satisfies
\begin{equation}
\mu^{\bullet}\,\int_0^{\infty} e^{- \beta t} \, dG^{\bullet}(t) = 1\:,\label{malthus}
\end{equation}
where $\mu^{\bullet}$ is the expected number of {\it mutant} first-generation offspring of a newly mutated bacterium. We assume that both $\mu^{\bullet}$ (we frequently write $\mu$ instead of $\mu^{\bullet}$) and $\mu^{\circ}$ are finite, that neither $G^{\circ}(\cdot)$ nor $G^{\bullet}(\cdot)$ are lattice-like, and that $G^{\circ}(0) = G^{\bullet}(0) = 0$. Then at each given moment in time, there are only finitely many bacteria (or \lq particles') alive \cite{sch}, and provided that the population has not become extinct (which is the case with probability $1 - q > 0$, and which we may easily assume), the population grows essentially as $\sim e^{\beta t}$ with $t$. The {\it expected} number of mutations in the population grows as $\sim \rho e^{\beta t}$. One of the basic assumptions of fluctuation analysis is that, as one allows the population to grow for longer and longer periods of time, the mutation rate in the population has always been so small that
\[
\rho e^{\beta t} \to m
\]
as $t\to\infty$, where $m < \infty$. This is the parameter entering into Equation \reff{straight}, and the one of predominate interest to fluctuation analysts.

The problem is thus to extract $m$ from experimental data under the model assumption \reff{straight}. To this end, we impose the following restriction on $f(s)$: Define
\begin{equation}
h(s) := \frac{1 - f(s)}{1 - s} \:, \label{h}
\end{equation}
where $f(s)$ is the PGF in \reff{f}. Our standing assumption is that, for $s$ in the neighbourhood of zero,
\begin{flushright}
$\mu^{\bullet} - h(1 - s) = s^{\omega}(-\log s)^{-\alpha} \cL(- \log s) \:,$ \hspace{2.2cm} $(\star)$
\end{flushright}
where either $\omega > 0$ and $\alpha\in(-\infty,\infty)$, or $\omega = 0$ and $\alpha \geq 2$. $\cL(\cdot)$ is slowly varying at infinity, and in case $\omega = 0$ and $\alpha = 2$ must be such that the integral
\begin{equation}
\int_1^{\infty} \frac{\cL(t)}{t}\,dt\label{cl}
\end{equation}
is finite. $(\star)$ is slightly stronger than the \lq $X \log X$'-condition of the Kesten-Stigum theorem. Via Tauberian arguments, it translates into $\sum_{k = n}^{\infty} k p_k \sim n^{\omega-1}(\log n)^{- \alpha} \cL(\log n)$ for $\omega<1$, else there is little intuitive about it. It is inspired by a paper of Uchiyama \cite{u}.

We shall show that if $(\star)$ holds and if the mutation under consideration is neutral in the sense that the Malthusian parameter $\beta$ is determined by {\it either} Equation \reff{nmmalthus} or \reff{malthus}, the law in \reff{basic} is in the normal domain of attraction of a stable distribution with index 1. Now the model considered by Lea and Coulson \cite{lc} in their celebrated formulation of the LDD {\it is} one of a neutral mutation and meets $(\star)$ with $h(s) = 1 + s$. Therefore its law is attracted to a stable distribution with index 1, and any part of the statistical machinery developed in Lea and Coulson's paper which rests on this fact alone can be can be applied no matter how cell proliferation proceeds in detail. The catch is that one will, in general, not obtain an estimate for the mutation rate itself, but only for the mutation rate times a factor
\begin{equation}
n^{\bullet}_1 := \frac{\mu - 1}{\beta \mu^2 \int_0^{\infty} t e^{- \beta t} dG^{\bullet}(t)} \:. \label{n1}
\end{equation}
It is seen that $n^{\bullet}_1$ is completely determined by the life-time distribution $G^{\bullet}(t)$ and the expected number of offspring ($\mu =\mu^{\bullet}$) of mutant bacteria. If either is unknown, $n_1^{\bullet}$ must be determined from auxiliary experiments. This may seem a drawback of the method advocated in this paper, but it is actually a problem inherent to fluctuation analysis itself: It seems intuitive that by mere inspection of the distribution of mutant colonies in an experiment, it might be difficult to disentangle how mutation and cell proliferation work together to produce the actual total number of mutants each population contains (think of a population where mutants are produced at a higher rate but it takes longer for mutant cells to divide). From this point of view, a simple proportionality between the estimated and the \lq actual' mutation rate is probably as transparent a relation as one could wish for.

We are now ready to state our main
\begin{theorem}\label{main}
Suppose that the PGF for the number of mutants $r$ in a large culture of bacteria is given as in \reff{basic}. Suppose furthermore that the life-time distributions of neither mutant nor non-mutant bacteria are lattice-like, and that the mutation under consideration is neutral in the sense that the Malthusian parameter for a population of exclusively non-mutant bacteria is the same as that for a population of exclusively mutant ones. Then, with $n_1^{\bullet}$ as given in Equation \reff{n1}, there exists $\delta\in(-\infty,\infty)$ such that, as the expected number $m$ of mutations in the population tends to infinity,
\begin{displaymath}
\xi_m := \frac{r}{n_1^{\bullet} m} - \log n_1^{\bullet} m - \frac{\delta}{n_1^{\bullet}}
\end{displaymath}
converges in distribution to a random variable $\xi$ with characteristic function $\ex(e^{i\xi\theta}) = (-i\theta)^{-i\theta}$ iff condition $(\star)$ is fulfilled.
\end{theorem}
We shall come to the proof in a moment, but let us remark now that the gist of the theorem lies in the scaling and centering of the distribution of mutants by $n_1^{\bullet}m$ and $\log n_1^{\bullet}m$, respectively, and in that the limiting distribution has been identified. This allows one to employ a maximum-likelihood method for estimating $n_1^{\bullet}m$ and $\delta/n_1^{\bullet}$ as described in \cite{wa2}. As of yet, the exact value of $\delta/n_1^{\bullet}$ does not seem to be of much importance, but the fact that, ultimately, the experiment will yield an estimate for the {\it product} of $n_1^{\bullet}$ times the mutation rate rather than the mutation rate itself, does. We defer further discussion of this topic to the final section of the paper, and proceed to the\\

\noindent{\it Proof of Theorem \ref{main}}. By Lemma \ref{thbh}, whose proof will be given in Section \ref{sbh} below,
\begin{equation}
m \big(g(s) - 1\big) =: m \big(1 - s\big)\big(n_1^{\bullet} \log(1 - s) + \delta(s)\big)\label{bigstep}
\end{equation}
for all $s \in [0, 1]$, where $\delta(s)$ is bounded and has a finite limit as $s\to 1$. Therefore we obtain for the characteristic function $\psi_m(\theta)$ of the random variable $\xi_m$,
\begin{eqnarray}
\lefteqn{\log \psi_m(\theta) = m \Big(1 - e^{i\theta(n_1^{\bullet} m)^{-1}} \Big) \times} \nonumber \\
&& \Big(n_1^{\bullet} \log\big(1 - e^{i\theta(n_1^{\bullet} m)^{-1}}\big) + \delta(e^{i\theta(n_1^{\bullet} m)^{-1}})\Big) - i\theta\log n_1^{\bullet}m - i\theta\,\frac{\delta}{n_1^{\bullet}} \nonumber\\
&=& - i\theta\log \big(-i\theta(n_1^{\bullet} m)^{-1} + O(m^{-2})\big) + i \theta\,\frac{\delta(e^{i\theta(n_1^{\bullet} m)^{-1}})}{n_1^{\bullet}}\nonumber\\
&&-i\theta\log n_1^{\bullet} m - i\theta\,\frac{\delta}{n_1^{\bullet}}  + O(m^{-2})\:, \nonumber
\end{eqnarray}
so that letting $m \to \infty$, we find $\psi_m(\theta) \to \psi(\theta) := (-i\theta)^{-i\theta}$ if we choose $\delta := \lim_{s\to 1}\delta(s)$. This already concludes the proof of the theorem.\hfill $\Box$\\

Of course, the largest part work with the proof will be to show that the claim concerning the continuity of $\delta(s)$ at $1^-$ made in \reff{bigstep} holds. The problem is to find a workable method of how to compare the function $g(s)-1$ with what is written on the right-hand side of \reff{bigstep}. As we shall see, there is one more or less natural way to do this if both mutant and non-mutant cells grow by way of an ordinary Galton-Watson process, but this is an example which does not meet the requirements of the theorem (the life-time distributions are very lattice-like in this instance), and we shall produce an example when a statement comparable to Theorem \ref{main} does not survive in the Galton-Watson world. One must therefore exercise some caution so not to prove what isn't true.
\section{The Galton-Watson Scenario} \label{prel}
Suppose the generating function $g(s)$ entering into Equation \reff{basic} were given by
\begin{equation}
g(s) :=: g_f(s) := \frac{\mu - 1}{\mu} \sum_{i=0}^{\infty}\frac{f_i(s)}{\mu^i}\:,\label{cg}
\end{equation}
where $f_i(\cdot)$ is for the $i$-th iterate of the generating function \reff{f}, and $f_0(s) := s$. We now prove the following
\begin{prop} \label{cthgw}
Let $f(s) = s^{\mu}$ for some integer $\mu \geq 2$, and let $g(s)$ be given as in Equation \reff{cg}. Define
\begin{displaymath}
\gamma(s) :=: \gamma_f(s) := \frac{1 - g(s)}{1 - s}\:.
\end{displaymath}
Then $\lim_{s \rightarrow 1^-} \gamma(s) + \kappa \log(1 - s)$ does not exist in $(-\infty,\infty)$ for any choice of $\kappa$.
\end{prop}
{\it Proof}. One easily checks that
\[
\mu g(s) = g\circ f(s) + (\mu-1)s
\]
and
\[
\mu \gamma(s) = \frac{\mu - g\circ f(s) - (\mu-1)s}{1-s} = \mu-1 + h(s)\,\gamma\circ f(s)\:.
\]
Let us now introduce a function $\vt(\cdot)$ and a change of variables $s =: e^{-e^{-t}}$ such that
\[
\vt(t) := \gamma(e^{-e^{-t}})\:.
\]
Thus, if $f(s) = s^\mu$, we find
\[
\vt(t) = \frac{\mu-1}{\mu} + \frac{\bar{h}(t)}{\mu}\,\vt(t-\log \mu)\:,
\]
where $\bar{h}(t) := h(e^{-e^{-t}})$. Thus we see that over a range of $\log\mu$, $\vt(t)$ increases by about $1-\mu^{-1}$, so it makes sense to ask whether $\vt(t)$ might not be approximately linear in $t$; or, to be more precise, to ask whether
\begin{equation}
\tau(t) := \vt(t) - \frac{\mu-1}{\mu\log\mu}\,t = \tau(t-\log\mu) - \frac{\mu - h(t)}{\mu}\,\gamma(t-\log \mu)
\end{equation}
has a limit as $t\to\infty$. Because $h(\cdot)$ meets $(\star)$ if $f(s) = s^{\mu}$, we have
\[
\mu - \bar{h}(t)\leq {\mu\choose 2} e^{-t}
\]
which implies that, indeed, $\lim_{t\to\infty}\tau(t)$ exists at least if $t$ runs through integer multiples of $\log\mu$. This means that the lemma has been proved except for $\kappa = \frac{\mu-1}{\mu\log\mu}$ and if
\[
\delta(s) := \gamma(s) + \frac{\mu-1}{\mu\log\mu}\,\log( - \log s)
\]
had a finite limit as $s\to 1^-$. We show that this is not the case: In fact, one easily checks that the coefficient of $s^k$ in $\gamma(s)$ equals $\mu^{-i}$ for $\mu^{i-1}\leq k \leq \mu^i - 1$, and since $\mu-1 \geq \log\mu$ for any non-negative $\mu$, it follows that
\begin{eqnarray}
\lefteqn{\sum_{k=\mu^{i-1} (\mu-1)/\log\mu}^{\mu^i-1} \frac{1}{\mu^i} - \frac{\mu-1}{\mu\log\mu}\frac{1}{k}}\nonumber\\
&\simeq& 1 - \frac{\mu-1}{\mu\log\mu}\left(1 - \log\frac{\mu-1}{\mu\log\mu}\right)\:,\nonumber
\end{eqnarray}
so that the sequence of partial sums of the first $n$ coefficients of $\delta(s)$ cannot converge. However, these coefficients being of the order $1/k$ for the $k$-th coefficient, they would have to if $\delta(s)$ converged for $s\to\infty$, by a Tauberian theorem of Hardy and Littlewood \cite{h}. But since
\[
\log(-\log s) - \log(1 - s) = \log\left(\frac{-\log s}{1-s}\right) \to 0
\]
as $s\to 1$, this contradiction concludes the proof of the lemma.\hfill$\Box$\\

In contrast to this, we have
\begin{prop} \label{frrat}
Let
\[
f(s) = \frac{s}{s + \mu(1 - s)}
\]
for some real number $\mu > 1$, and let $\gamma(s) :=: \gamma_f(s)$ be as in Proposition \ref{cthgw}. Then
\[
\lim_{s\to 1^-} = \gamma(s) + \frac{\mu-1}{\mu\log\mu}\,\log(1-s)
\]
exists in $(-\infty,\infty)$.
\end{prop}
{\it Proof}. First observe that
\[
\mu - h(1-s) = s\,\frac{(\mu-1)\mu}{1+(\mu-1)s}\:,
\]
so condition $(\star)$ is fulfilled in this case, and that $f(\cdot)$ is of the fractional rational type. Hence we have for its $i$-th iterate \cite{an},
\[
\frac{1 - f_i(s)}{\mu^i} = \frac{1-s}{\mu^i + (1-\mu^i)s}\:,
\]
which implies
\[
\gamma(s) = \frac{\mu-1}{\mu}\,\sum_{i=0}^{\infty} \frac{1}{\mu^i + (1-\mu^i)s} = \frac{\mu-1}{\mu}\,\sum_{k=0}^{\infty} s^k \sum_{i=0}^{\infty}\frac{1}{\mu^i}\left(1-\frac{1}{\mu^i} \right)^k\:.
\]
Consider now the sum of the first $n-1$ coefficients. Except for the factor $1- \mu^{-1}$, it is
\[
\sum_{k=0}^{n-1}\sum_{i=0}^{\infty} \frac{1}{\mu^i} \left(1 -\frac{1}{\mu^i}\right)^k = \sum_{i=0}^{\infty} \frac{1}{\mu^i} \frac{1 - (1 - \mu^{-i})^n}{1 - (1 - \mu^{-i})} = \sum_{i=0}^{\infty}1 - (1 - \mu^{-i})^n\:.
\]
By Euler's formula, this equals
\begin{eqnarray}
\lefteqn{\sum_{i=0}^{\infty} 1 - (1 - \mu^{-i})^n = \int_0^{\infty} \big(1 - (1-\mu^{-u})^n\big)\,du}\nonumber\\
&&-\:n \log \mu \int_0^{\infty} (u - [u]) \mu^{-u}(1-\mu^{-u})^{n-1}\,du + 1\:.\nonumber
\end{eqnarray}
We re-write the middle term as
\begin{eqnarray}
\lefteqn{\log \mu \int_0^{\infty} (u - [u]) n \mu^{-u}\left(1- \frac{n\mu^{-u}}{n}\right)^{n-1}\,du}\nonumber\\
&=& \log \mu \int_{-\log n/\log \mu}^{\infty} (v - [v]) \mu^{-v}\left(1- \frac{\mu^{-v}}{n}\right)^{n-1}\,dv\:,\nonumber
\end{eqnarray}
which readily follows after a change of variables $u =: v + \log n/\log \mu$. But this last expression tends to
\[
\log \mu \int_{-\infty}^{\infty} (v - [v]) \mu^{-v}\exp(-\mu^{-v})\,dv < \infty
\]
as $n\to\infty$. Now one easily checks by expanding the integrand that
\[
\int_0^{\infty} \big(1 - (1 - \mu^{-u})^n\big)\,du = \frac{H_n}{\log\mu}\:,
\]
where $H_n$ is the $n$-th harmonic number. It follows that the discounted partial sums
\[
\sum_{k=0}^{n-1}\sum_{i=0}^{\infty} \frac{1}{\mu^i} \left(1 -\frac{1}{\mu^i}\right)^k - \frac{H_n}{\log\mu}
\]
converge, which by Abel's theorem already implies the lemma.\hfill$\Box$\\

This is something of a mess, because, as we have seen during the proof of Theorem \ref{main}, the theorem actually asserts that the function $\gamma(s)$ behaves essentially as $-n_1^{\bullet}\log(1-s)$ in the neighbourhood of 1, and by now we have both an example where it does, as well as an expample where it doesn't. The next step is therefore to suppose that $g(s)$ is again as given in Equation \reff{cg}, but with $f(\cdot)$ a PGF which isn't either $s^{\mu}$ or $\frac{s}{s + \mu(1 - s)}$. Then, with $\gamma(s)$ as defined in Proposition \ref{cthgw}, and with $\kappa = \frac{\mu-1}{\mu\log\mu}$, one finds that
\begin{eqnarray}
\lefteqn{\gamma(s) + \kappa \log(1 - s) = \frac{\mu - 1}{\mu} \sum_{i=0}^{\infty} \frac{1}{\mu^i} \, \frac{1 - f_i(s)}{1 - s}} \label{triple} \\
&&- \: \frac{\kappa}{\mu} \, \sum_{i=0}^{\infty} \frac{h \circ f_i(s)}{\mu^i} \, \frac{1 - f_i(s)}{1 - s} \, \log \big(h \circ f_i(s)\big) \nonumber \\
&&+ \: \frac{\kappa}{\mu} \, \sum_{i=0}^{\infty} \frac{\mu - h \circ f_i(s)}{\mu^i} \,  \frac{1 - f_i(s)}{1 - s} \, \log\big(1 - f_i(s)\big) \:. \nonumber
\end{eqnarray}
The way to see this is to use the definition \reff{h} of $h(\cdot)$ and to check that the second and third sum, taken together, telescope. Now fix $s \in (q, 1)$, and denote by $f_{-1}(s)$ the inverse of $f(s)$. Then $f_{-n}(s) := f_{-1} \circ f_{1-n}(s) \to 1$ as $n \to \infty$ (because $f'(1) = \mu > 1$), and one deduces from \reff{triple} that
\begin{eqnarray}
\lefteqn{\gamma \circ f_{-n}(s) + \kappa \log\big(1 - f_{-n}(s)\big) = \frac{\gamma \circ f(s) + \kappa \log\big(1 - f(s)\big)}{\mu^n \big(1 - f_{-n}(s)\big)}} \label{tobv} \\
&+& \frac{1}{\mu^n \big(1 - f_{-n}(s)\big)} \bigg( \frac{\mu - 1}{\mu} \sum_{i=1}^n \mu^i \big(1 - f_{-i}(s)\big) \nonumber \\
&&- \: \frac{\kappa}{\mu} \, \sum_{i=1}^n \mu^i \big(1 - f_{-i}(s)\big) h \circ f_{-i}(s) \log \big(h \circ f_{-i}(s)\big) \nonumber \\
&&+ \: \frac{\kappa}{\mu} \, \sum_{i=1}^n \mu^i \big(1 - f_{-i}(s)\big)\big(\mu - h \circ f_{-i}(s)\big) \log\big(1 - f_{-i} (s)\big) \bigg) \:. \nonumber
\end{eqnarray}
It can be shown \cite{an} that under condition $(\star)$, 
\begin{equation}
\mu^n \big(1 - f_{-n}(s)\big) \to Q^{\star}(s) \label{qs}
\end{equation}
as $n \to \infty$, with $0 < Q^{\star}(s) < \infty$ for all $s \in (q,1)$. Therefore, since
\[
\mu - h \circ f_{-i}(s)\sim \mu - h\big(1 - \mu^{-i} Q^{\star}(s)\big) \sim \frac{1}{\mu^{i\omega}} i^{-\alpha}\cL(i)\:,
\]
the final right-hand term in the above equation is of the order
\begin{displaymath}
\sim\sum_{i=1}^n \, \frac{\cL(i)}{\mu^{i\omega}i^{\alpha - 1}}\:,
\end{displaymath}
and because $(s - 1) \log(1 - s) = s + O(s^2)$ for $s$ in the neighbourhood of zero, we also have
\begin{eqnarray}
\lefteqn{\mu\,\sum_{i=1}^n \Big(\log \mu - \log \big(h \circ f_{-i}(s)\big)\Big) \: + \: \sum_{i=1}^n \log h \circ f_{-i}(s) \big(\mu - h \circ f_{-i}(s)\big)}\nonumber \\
&\sim&\sum_{i=1}^n \Big(- \log \big(1 - f_{-i}(s)\big)\Big)^{- \alpha}\cL\Big(- \log \big(1 - f_{-i}(s)\big)\Big)\sim \sum_{i=1}^n\frac{\cL(i)}{\mu^{i\omega}i^{\alpha}}\:.\nonumber
\end{eqnarray}
Overall, we find that $\lim_{n \to \infty} \gamma \circ f_{-n}(s) + \kappa \log\big(1 - f_{-n}(s)\big)$ is of the order $\sim\sum_{i=1}^n \cL(i) \mu^{-i\omega} i^{1 - \alpha}$, which is bounded because of $(\star)$. The problem is that this only implies the existence of
\[
\lim_{n\to\infty}\gamma\circ f_{-n}(s) + \frac{\mu-1}{\mu\log\mu}\log\big(1 - f_{-n}(s)\big)
\]
for every fixed value of $s\in(q,1)$, but doesn't say anything about whether this limit is independent of $s$ or not. It turns out, however, that for Bellman-Harris processes with a non-lattice life-time distribution of particles, one can turn the above argument into an honest proof of Theorem \ref{main}. We will treat the example which is closest at hand in the following section.
\section{Markov Branching in Continuous Time}
We will now show that the argument given towards the end of the previous section goes through for continuous-time Markov models of clonal expansion. The construction goes as follows: Define $\delta' :=: \delta/\beta$ for some $\delta > 0$, and let
\begin{equation}
\phi(s) :=: \phi(s,\delta) := \frac{1}{\delta'}\int_0^{\delta'} e^{- \beta u} F_u(s)\,du\:,
\end{equation}
and
\begin{equation}
\varphi(s) :=: \varphi(s,\delta) := \frac{\phi(1) - \phi \circ F_{\delta'}(s)}{\phi(1) - \phi(s)} \:.
\end{equation}
Except for a norming factor, $\phi(s)$ is a generating function, and we will find that $\varphi(s)$ corresponds to the function $h(s)$ defined in Equation \reff{h} in a way which is sufficient for our purposes. Also,
\[
\phi(1) = \frac{1}{\beta\delta'}\,(1 - e^{-\beta\delta'}) = \frac{1}{\delta}\,(1 - e^{-\delta})\:,
\]
$\phi'(1) = 1$, and 
\begin{displaymath}
\varphi(1) = \lim_{s \to 1} \frac{1}{\delta'} \int_0^{\delta'} e^{- \beta u} \frac{1 - F_{u + \delta'}(s)}{1 - s} = \frac{1}{\delta'} \int_0^{\delta'} e^{- \beta u} e^{\beta (u + \delta')} \, du = e^{\beta\delta'} = e^{\delta}\:,
\end{displaymath}
{\it provided that the Malthusian parameters for mutant and non-mutant cells are the same}. This is where the assumption of no selection on mutant cells enters the calculations. Observe also that
\[
\frac{1 - F_{u+\delta'}(s)}{1 - F_u(s)} = 
\frac{1 - F_{\delta'}\circ F_u(s)}{1 - F_u(s)} < e^{\delta}\:,
\]
so that $e^{\delta} - \varphi(s)\geq 0$. Now, since the generating functions $F_u(s)$ have the semi-group property, we find
\begin{eqnarray}
\beta \int_0^{\infty} e^{- \beta u} F_u(s) \, du &=& \beta \int_0^{\delta'} e^{- \beta u} F_u(s) \, du \:+\: \beta \int_{\delta'}^{\infty} e^{- \beta u} F_u(s) \, du
\nonumber \\
&=& \beta\delta'\phi(s) \:+\: \frac{\beta}{e^{\beta\delta'}} \int_0^{\infty} e^{- \beta u} F_u \circ F_{\delta'}(s) \, du \:, \nonumber
\end{eqnarray}
so that
\begin{equation}
g(s) = \delta\sum_{i=0}^{\infty} \frac{\phi \circ F_{i \ast \delta'}(s)}{e^{i\delta}}\:,
\end{equation}
where we have denoted by $F_{i \ast \delta}(s)$ the $i$-th iterate of the PGF $F_{\delta}(s)$. Note that
\[
g(1) = \delta\sum_{i=0}^{\infty} \frac{\phi(1)}{e^{i\delta}} = \frac{\delta\phi(1)}{1 - e^{-\delta}} = 1\:,
\]
as it must be. Also,
\[
\gamma(s) = \frac{1-g(s)}{1-s} = \delta\sum_{i=0}^{\infty} \frac{1}{e^{i\delta}}\,\frac{\phi(1) - \phi \circ F_{i \ast \delta'}(s)}{1-s}\:, 
\]
and therefore
\begin{eqnarray}
\lefteqn{\gamma(s) + \frac{\phi(1)-\phi(s)}{1-s}\log\big(\phi(1) - \phi(s)\big) = \delta \sum_{i=0}^{\infty} \frac{1}{e^{i\delta}} \, \frac{\phi(1) - \phi \circ F_{i \ast \delta'}(s)}{1 - s}}
\label{triplem} \\
&&- \: \frac{1}{e^{\delta}}\sum_{i=0}^{\infty} \frac{\varphi \circ F_{i \ast \delta'}(s)}{e^{i\delta}} \, \frac{\phi(1) - \phi \circ F_{i \ast \delta'}(s)}{1 - s} \, \log \big(\varphi \circ F_{i \ast \delta'}(s)\big) \nonumber \\
&&+\:\frac{1}{e^{\delta}}\sum_{i=0}^{\infty} \frac{e^{\delta} - \varphi \circ F_{i \ast \delta'}(s)}{e^{i\delta}} \,  \frac{\phi(1) - \phi \circ F_{i \ast \delta'}(s)}{1 - s} \, \log\big(\phi(1) - \phi \circ F_{i \ast \delta'}(s)\big) \:. \nonumber
\end{eqnarray}
Again, one checks this via inserting for $\log\varphi(\cdot)$ into the second sum according to the definition of $\varphi(\cdot)$ and rearranging terms (which causes no trouble because of the absolute convergence of the sums), whereupon one ends up with a single sum which telescopes to
\[
\frac{\phi(1) - \phi \circ F_0(s)}{1 - s} \, \log\big(\phi(1) - \phi \circ F_0(s)\big) = \frac{\phi(1)-\phi(s)}{1-s}\log\big(\phi(1) - \phi(s)\big)\:.
\]
The equality \reff{triplem} established, it follows easily from there that, if we denote by $F_{-i\ast\delta'}(\cdot)$ the inverse of the function $F_{i\ast\delta'}(\cdot)$,
\begin{eqnarray}
\lefteqn{\gamma\circ F_{-n\ast\delta'}(s) + \frac{\phi(1)-\phi\circ F_{-n\ast\delta'}(s)}{1-F_{-n\ast\delta'}(s)} \log\big(\phi(1) - \phi\circ F_{-n\ast\delta'}(s)\big)}\nonumber
\\
&=&\frac{1-s}{e^{n\delta}\big(1 - F_{-n\ast\delta'}(s)\big)}\Big(\gamma(s) + \log\big(\phi(1) - \phi(s)\big)\Big)\nonumber\\
&&+\:\frac{1}{e^{n\delta}\big(1 - F_{-n\ast\delta'}(s)\big)}\sum_{i=1}^n e^{i\delta}\,\big(\phi(1) - \phi \circ F_{-i \ast \delta'}(s)\big)\Delta_i(s)\:, \label{tokin}
\end{eqnarray}
where
\begin{eqnarray}
\lefteqn{\Delta_i(s) = \delta - \frac{\varphi \circ F_{-i\ast \delta'}(s)}{e^{\delta}}\,\log \big(\varphi \circ F_{-i\ast \delta'}(s)\big)}\nonumber
\\
&&+\:\frac{e^{\delta} - \varphi \circ F_{-i\ast \delta'}(s)}{e^{\delta}}\,\log\big(\phi(1) - \phi \circ F_{-i\ast \delta'}(s)\big) \:.\label{ds}
\end{eqnarray}
Thus, if we could show that the limit $n\to\infty$ in \reff{tokin} exists and is independent of $\delta$ and $s$, we would have made a big step towards a proof of the following
\begin{lemma} \label{thctm} Under condition $(\star)$, and provided that the life-time distribution $G^{\bullet}(t)$ of mutant cells is exponential (that is, $G^{\bullet}(t) = 1 - e^{- \lambda t}$ for some $\lambda > 0$) we have that
\begin{displaymath}
\lim_{s \to 1^-} \gamma(s) + \log(1 - s)
\end{displaymath}
exists and is bounded away from $\pm \infty$.
\end{lemma}
{\it Proof}. We first prove that the limit $n\to\infty$ in \reff{tokin} exists. Because
\[
e^{i\delta}\,\big(\phi(1) - \phi \circ F_{-i \ast \delta'}(s)\big) = e^{i\delta}\big(1 - F_{-i\ast \delta'}(s)\big)\,\frac{\phi(1) - \phi \circ F_{-i\ast \delta'}(s)}{1 - F_{-i\ast \delta'}(s)}
\]
converges as $i\to\infty$ (see \cite{an}, and observe that $\phi'(1)=1$), what we need to show is that the $\Delta_i(s)$ tend to zero sufficiently quickly. To do this, write $M_t$ for the expected size of a mutant clone of age $t$, and set $H_t(s) := \frac{1 - F_t(s)}{1 - s}$. Then $M_t = H_t(1)$, and
\begin{eqnarray}
\lefteqn{M_t - H_t(s)} \nonumber \\
&=& \lambda \mu \int_0^t M_{t-u} \, e^{- \lambda u} \,du - \lambda \int_0^t h \circ F_{t - u}(s) H_{t-u}(s) \, e^{- \lambda u)} \,du \nonumber \\
&=& \lambda \int_0^t \big(\mu - h \circ F_{t-u}(s)\big) H_{t-u}(s) \, e^{- \lambda u} \,du \nonumber \\
&&+ \: \lambda \mu \int_0^t \big(M_{t-u} - H_{t-u}(s) \big) \, e^{- \lambda u} \,du \:, \nonumber
\end{eqnarray}
which can be solved by means of Laplace transforms:
\begin{displaymath}
M_t - H_t(s) = \lambda \int_0^t \big(\mu - h \circ F_{t-u}(s)\big) H_{t-u}(s) \, e^{\beta u} \, du \:.
\end{displaymath}
Then, since $M_t = e^{\beta t}$, it is straightforward to check that
\begin{eqnarray}
\lefteqn{0\leq e^{\delta} - \varphi(s)\leq\frac{1}{\delta'} \int_0^{\delta'}e^{-\beta u} \big(M_{u + \delta'} - H_{u+\delta'}(s)\big)\,du} \nonumber \\
&=& \frac{\lambda}{\delta'}\int_0^{\delta'}e^{-\beta u}\int_0^{u+\delta'} \big(\mu - h \circ F_{u+\delta'-v}(s)\big) H_{u+\delta'-v}(s) \, e^{\beta v}\,dv\,du\nonumber\\
&\leq& \lambda(\mu - h \circ F_{2\delta'}(s)\big)\:.\nonumber
\end{eqnarray}
But it follows easily from Jensen's inequality that $F_t(s) \geq s^{M_t} = s^{e^{\beta t}}$, which by $(\star)$ immediately gives
\[
e^{\delta} - \varphi(1-s) \leq \lambda e^{2\delta\omega} s^{\omega} \big(-\log s-2\delta\big)^{-\alpha} \cL\big(-\log s - 2\delta \big)\:.
\]
This is good enough for us: because
\[
\log \big(\varphi \circ F_{-i\ast \delta'}(s)\big) = \delta + \log \left(1 - \frac{e^{\delta} - \varphi \circ F_{-i\ast \delta'}(s)}{e^{\delta}}\right)
\]
differs, in the neighbourhood of 1, from $\delta$ only by an amount of order $e^{\delta} - \varphi \circ F_{-i\ast \delta'}(s)$, the $\Delta_i(s)$ as given by \reff{ds} tend to zero at least as quickly as $\sim\cL(i)i^{1-\alpha}$, which proves that the sum in Equation \reff{tokin} converges. But
\[
1 - \frac{\phi(1)-\phi(s)}{1-s} = \frac{1}{\delta'} \int_0^{\delta'} e^{- \beta u} \big(M_u - H_u(s)\big)\,du\:,
\]
which by the above implies that
\[
\frac{\phi(1)-\phi(s)}{1-s} \log(1-s) - \log(1-s)\to 0
\]
as $s\to 1$. All in all, we know now that
\begin{equation}
\lim_{n\to\infty}\gamma\circ F_{-n\ast\delta'}(s) + \log\big(1 - F_{-n\ast\delta'}(s)\big)\label{ctmking}
\end{equation}
converges for every $\delta'>0$ and every $s\in(q,1)$. But $F_t(\cdot)$ is continuous in $t$; and since $\gamma(s)$ and $\log(1-s)$ are continuous for $s\in[0,1)$, an application of Kingman's lemma \cite{king} shows that the limit in \reff{ctmking} is independent of $\delta$. Since we can find for any $s\leq u$, $s,u\in (q,1)$ some suitable $n$ and $\delta$ such that $s = F_{n\ast\delta'}(u)$, this limit is also independent of $s$, which concludes the proof of Lemma \ref{thctm}. \hfill $\Box$
\subsection{A Second Proof}
For completeness, we now give a second proof of Lemma \ref{thctm} which makes use of the forward picture for the continuous time Markov process, namely \cite{an}
\[
\frac{\partial}{\partial t} F_t(s) = u(s) \frac{\partial}{\partial s}F_t(s)\:,
\]
where
\[
u(s) = \lambda\big(f(s)-s\big)\:.
\]
Then we have by definition of $g(s)$,
\begin{eqnarray}
\lefteqn{g(s) = -e^{-\beta u} F_u(s)\Big\vert_{u=0}^{\infty}+\int_0^{\infty}e^{-\beta u} \frac{\partial}{\partial u} F_u(s)\,du}\nonumber\\
&=& s -  \lambda\big(s - f(s)\big)\int_0^{\infty}e^{-\beta u} \frac{\partial}{\partial u}F_u(s)\,du = s - \frac{\lambda}{\beta}\big(s - f(s)\big)g'(s)\:,\nonumber
\end{eqnarray}
which, as one readily checks, implies
\begin{equation}
\gamma(s) = 1 + \bar{z}(s)\big(\gamma(s) - \gamma'(s)(1-s)\big)\label{diffg}
\end{equation}
where
\begin{equation}
\bar{z}(s) := \frac{\lambda}{\beta}\frac{s- f(s)}{1-s} = \frac{h(s) - 1}{\mu-1}\:.\label{z}
\end{equation}
Now define $t := -\log(1-s)$ and a function $\vt(\cdot)$ of $t$ such that $\gamma(s) =: \vt\big(-\log(1-s)\big)$. Then
\[
\vt(t) = 1 + z(t)\big(\vt(t) - \vt'(t)\big)\:,
\]
and $z(t) := \bar{z}(1 - e^{-t})$. This can be solved by standard methods. We have
\[
\vt(t) = \vt_0(t) + \vt_0(t) \int_{\vs}^t \frac{1}{\vt_0(u)z(u)}\,du\nonumber
\]
where
\[
\vt_0(t) := \vt(\vs)\exp\left(\int_{\vs}^s\frac{z(u)-1}{z(u)}\,du\right)
\]
is the solution of the homogeneous equation corresponding to \reff{diffg}, and $\vs \in (q,1)$. Now write
\begin{eqnarray}
\vt(t) - t + \vs &=& \vt_0(t) + \int_{\vs}^t \frac{\vt_0(t) - z(u)\vt_0(u)}{\vt_0(u)z(u)}\,du\nonumber\\
&=& \vt_0(t) + \int_{\vs}^t \frac{\vt_0(t) - \vt_0(u)}{\vt_0(u)z(u)}\,du - \int_{\vs}^t \frac{z(u)-1}{z(u)}\,du\:,\nonumber
\end{eqnarray}
and observe that $\vt_0(t) - \vt_0(u)$ is at worst of order
\[
0 \geq \vt_0(t) - \vt_0(u)\geq \int_u^t\frac{z(v)-1}{z(v)}\,dv
\]
in magnitude, so
\[
0 \geq \int_{\vs}^t \frac{\vt_0(t) - \vt_0(u)}{\vt_0(u)z(u)} \,du \geq \frac{1}{\vt_0(\vs)z(\vs)} \int_{\vs}^t (u -\vs)\frac{z(u)-1}{z(u)}\,du\:,
\]
which by definition \reff{z} of $z(\cdot)$ and because of ($\star$) is of order
\[
 - \int_{\vs}^t (u - \vs)\frac{\mu - h(1 - e^{-u})}{h(1 - e^{-u})-1}\,du \sim - \int_{\vs}^t e^{-\omega u} \cL(u) u^{\alpha-1}\,du > - \infty\:.
\]
This concludes the proof, since now it follows that
\[
\gamma(s) = \vt\big(-\log(1-s)\big) = -\log(1-s) + \delta\big(-\log(1-s)\big)\:,
\]
and $\delta(\cdot)$ approaches some finite value as its argument tends to infinity. \hfill $\Box$\\

By now, we have completed the proof of Theorem \ref{main} in case the life-time distribution of mutant cells is exponential. From a practical point of view, this is what we wanted the least, since in a real-life situation, newborn cells will have to pass through a whole cell cycle before they can divide anew. Therefore, in the following section, we shall imitate the first proof of Lemma \ref{thctm} to obtain a comparable statement for general models of cell proliferation.
\section{The General Case} \label{sbh}
One difficulty with non-exponential life-time distributions is that, if one only looks at the number of particles $Z_t$ at each given instant, the process $\{Z_t\}_{t\geq 0}$ is not Markovian anymore. This implies that the generating functions $F_t(s) = \ex(s^{Z_t})$ will not form a semi-group, so that in particular, $F_{-u+t} \circ F_{-t}(s) \neq F_{-u}(s)$. Nevertheless, $\lim_{t \to \infty} F_{-u+t} \circ F_{-t}(s)$ exists and can be calculated as follows: Fix $s \in (q, 1)$, where $q$, again, is the unique fixed point of $f(s)$ in $[0, 1)$. It can be shown \cite{sch} that the product of $- \log F_{-t}(s)$ times the number of bacteria $Z_t$ that by time $t$ have originated from a bacterium born at time 0 converges almost surely to a non-degenerate random variable $Z$. Then, with
\begin{equation}
R(x) :=: R_s(x) := \lim_{t \to \infty} F_t \circ F_{-t}(s)^x \label{r}
\end{equation}
the Laplace transform of the random variable $Z$, we have \cite{sch}
\begin{equation}
\lim_{t \to \infty} F_{-u+t} \circ F_{-t}(s)^x = R(xe^{- \beta u}) = \int_0^{\infty} f \circ R(xe^{- \beta u}) \, dG(u) \:. \label{was}
\end{equation}
Consider now the number of bacteria $Z_t[x; y]$ aged less than $y$ at time $t$ that have originated from a single bacterium aged $x$ at time 0. Then
\begin{equation}
Z_t[x] := \int_0^{\infty} Z_t[x; dy] \label{zulu}
\end{equation}
is for the total number of progeny an individual aged $x$ has produced by time $t$, and it is easy to see (say, by considering the corresponding Laplace transforms) that $- \log F_{-t}(s)$ times this number converges in distribution to random variable $Z[x]$. Now for $\fr(x)$ a measurable function on $[0, \infty)$, define
\begin{displaymath}
\int_0^{\infty} \fr(y) Z_t[x; dy] := \sum_{i=1}^{Z_t[x]} \fr(y_i) \:,
\end{displaymath}
where the $y_i$ are the ages of individual bacteria in the population at time $t$, and set
\begin{equation}
\Phi_t[x](\fr) := \ex\bigg(\exp \Big({\textstyle\int_0^{\infty}\log\fr(y) Z_t[x; dy]}\Big)\bigg) \:. \label{gf}
\end{equation}
This is the definition of the generating functional (GF). The $\Phi_t[\cdot](\cdot)$'s {\it do} have the semi-group property in the sense that $\Phi_{t + r}[x](\fr) = \Phi_t[x]\big(\Phi_r[\cdot](\fr)\big)$. Furthermore, Schuh \cite{sch} has shown that there exist measurable functions $\fs_t(x)$ on $[0, \infty)$ such that
\begin{equation}
\Phi_r[x](\fs_t) = \fs_{t-r}(x) \label{recur}
\end{equation}
for all $0 \leq r \leq t < \infty$. The $\fs_t(x)$'s can be specified explicitly \cite{sch}: With
\begin{equation}
G_y(t) := \frac{G(y+t) - G(y)}{1 - G(y)}
\end{equation}
the residual life-time distribution function of a mutant bacterium aged $y$ (we write $G(t)$ instead of $G^{\bullet}(t)$ for brevity),
\begin{equation}
\fs_t(x) := \int_0^{\infty} f \circ R(e^{- \beta(t + u)}) \, dG_x(u) \:. \label{fs}
\end{equation}
Since $R(0) = 1$, the $\fs_t(x)$ tend to 1 as $t \to \infty$. We need an estimate for the rate of convergence. {\it From now on, we will always assume that $G(\cdot)$ is non-lattice}. We have the following
\begin{lemma} \label{straight} For $\fs_t(x)$ as given by Equation \reff{fs}, $R(\cdot)$ as given in Equation \reff{r}, and
\begin{displaymath}
V(x) := \int_0^{\infty} e^{- \beta u}\,dG_x(u)
\end{displaymath}
the reproductive value, we have
\begin{displaymath}
\mathcal{C}\left\vert\mu V(x) + e^{-\beta u}\frac{1 - \fs_{t-u}(x)}{\log R(e^{-\beta t})}\right\vert \leq \mu -  h\circ R(e^{-\beta (t-u)})
\end{displaymath}
for some constant $\mathcal{C}> 0$ and any $u \leq t$.
\end{lemma}
The proof of this lemma is deferred to the appendix. Observe that convergence is uniform in $x$. Next, we prove
\begin{lemma} \label{statt} For
\begin{displaymath}
\syG[x](\fr) := \beta \int_0^{\infty} e^{-\beta u} \Phi_u[x](\fr) \, du \:,
\end{displaymath}
with $\Phi_u[x](r)$ as given in Equation \reff{gf}, and $\bar{\fs}_t(x) := \exp\big(n_1^{\bullet} \log \fs_t(x)\big)$, we have that
\begin{displaymath}
\lim_{t \to \infty} -\:\frac{\syG[x]\big(R(e^{- \beta t})\big) - \syG[x](\bar{\fs}_t)}{\log R(e^{- \beta t})} = \beta \int_0^{\infty}\Big(e^{- \beta u} M_u[x] - n_1^{\bullet}\mu V(x)\Big)\,du\:.
\end{displaymath}
In particular, the limit exists in $(-\infty,\infty)$.
\end{lemma}
{\it Proof}. Note first that the lemma is trivial if $Z_u[x]$ explodes in finite time. We may therefore assume that $Z_u[x] < \infty$ for all $t\in[0,\infty)$. Also
\[
\lim_{t\to\infty} Z_t^{-1} \int_0^y Z_t[0; dy] \to A(y)
\]
almost surely, where $A(\cdot)$ is the limiting age distribution of particles \cite{sch}. Granted this fact, the same is obviously true for $Z_t[x]^{-1} \int_0^y Z_t[x; dy]$. Now for the difference $u-t$ fixed, we have
\begin{eqnarray}
\lefteqn{\int_0^{\infty}\log\bar{\fs}_t(y) Z_u[x; dy]}\nonumber\\
&=& n_1^{\bullet}\log R(e^{-\beta u})Z_u[x]\int_0^{\infty}\frac{ \log \fs_t(y) }{\log R(e^{-\beta t})}\,\frac{\log R(e^{-\beta t})}{\log R(e^{-\beta u})}\,\frac{Z_u[x; dy]}{Z_u[x]}\nonumber\\
&\to&-\:n_1^{\bullet 2} \mu e^{\beta(u - t)} Z[x]\int_0^{\infty} V(y)\,dA(y) = -\:n_1^{\bullet}e^{\beta(u - t)} Z[x]\nonumber
\end{eqnarray}
at least in distribution as $t\to\infty$, and one finds
\begin{eqnarray}
\lefteqn{\lim_{t\to\infty}\int_t^{\infty} \frac{e^{-\beta u} \Phi_u[x](\bar{\fs}_t)}{\log R(e^{-\beta t})}\,du}\label{chuck}\\
&=& -\int_0^{\infty} e^{-\beta u}\ex(e^{- n_1^{\bullet}e^{\beta u} Z[x]})\,du = \lim_{t\to\infty}\int_t^{\infty} \frac{e^{-\beta u} \Phi_u[x]\big(R(e^{-\beta t})\big)}{\log R(e^{-\beta t})}\,du\:.\nonumber
\end{eqnarray}
Now $e^{-x}$ is a completely monotonous function in $x$; thus $e^{-y}(y-x) \leq e^{-x} - e^{-y} \leq e^{-x}(y-x)$ for arbitrary $x$ and $y$. Because of this, and because of Equation \reff{chuck}, we can say that
\begin{eqnarray}
\lefteqn {\lim_{t\to\infty}-\:\frac{\syG[x]\big(R(e^{- \beta t})\big) - \syG[x](\bar{\fs}_t)}{\log R(e^{- \beta t})}} \nonumber \\
&=&\lim_{t\to\infty}\beta\int_0^t e^{- \beta u} \, \ex \Bigg(\mathcal{E}_t(u) Z_u[x] + n_1^{\bullet} \mathcal{E}_t(u) \int_0^{\infty} \frac{\log\fs_t(y)}{\log R(e^{- \beta t})} \, Z_u[x; dy]\Bigg) \, du\:, \nonumber
\end{eqnarray}
where $Z_u[x]$ is as in Equation \reff{zulu}. Now $\mathcal{E}_t(u)$ is of order $\exp(- e^{-\beta t} Z_u[x])$. Thus, $\mathcal{E}_t(u)$ tends to 1 as $t\to\infty$ for every fixed value of $u$, and the lemma will be proved once we have shown that
\[
\beta\int_0^t e^{- \beta u} \, \ex \Bigg(Z_u[x] + n_1^{\bullet} \int_0^{\infty} \frac{\log\fs_t(y)}{\log R(e^{- \beta t})} \, Z_u[x; dy]\Bigg) \, du
\]
converges as $t\to\infty$. By linearity of taking expectations, the above can be re-written as
\begin{eqnarray}
\lefteqn{\lim_{t\to\infty}\beta \int_0^t \Big(e^{- \beta u} M_u[x] - n_1^{\bullet} \mu V(x)\Big) \, du} \nonumber \\
&+& \beta n_1^{\bullet} \int_{u=0}^t \int_{y=0}^{\infty} e^{- \beta u} \bigg(\mu V(y) + \frac{\log\fs_t(y)}{\log R(e^{- \beta t})} \bigg) \, M_u[x; dy] \, du\:, \nonumber
\end{eqnarray}
since $V(y)$ is a right eigenfunction to the kernel $M_t[x; y] := \ex \big( Z_t[x; y] \big)$ with eigenvalue $e^{\beta t}$ \cite{sch}. But
\begin{eqnarray}
\lefteqn{\int_{u=0}^t \int_{y=0}^{\infty} e^{- \beta u} \bigg(\mu V(y) + \frac{\log\fs_t(y)}{\log R(e^{- \beta t})} \bigg)\,M_u[x; dy]\,du}\nonumber\\
&\sim& e^{-\beta \omega t} t^{-\alpha} \cL(t)\int_0^t e^{-\beta u} M_u[x]\,du\hspace{2.2cm}\nonumber
\end{eqnarray}
because of Lemma \ref{straight}, and this tends to zero for large values of $t$. It therefore remains to show that $e^{- \beta u} M_u[x] - n_1^{\bullet} \mu V(x)$ is integrable. We have \cite{sch}
\begin{eqnarray}
\lefteqn{e^{-\beta t} M_t[y] - n_1^{\bullet} \mu V(y) = e^{-\beta t}\big(1 - G_y(t)\big)}\nonumber\\
&&-\:n_1^{\bullet}\mu\int_t^{\infty} e^{-\beta u}\,dG_y(u)+\mu \int_0^t \big(e^{-\beta (t-u)} M_{t-u} - n_1^{\bullet})\,e^{-\beta u}\,dG_y(u)\:,\nonumber
\end{eqnarray}
where $M_t$ is short for $M_t[0]$, and one readily checks that
\[
\int_0^{\infty} e^{-\beta t}\big(1 - G_y(t)\big)\,dt = \frac{1 - V(y)}{\beta}\:,
\]
and
\[
\int_{t=0}^{\infty} \int_{u=t}^{\infty} e^{-\beta u}\,dG_y(u)\,dt = \int_0^{\infty} u e^{-\beta u}\,dG_y(u) < \infty\:.
\]
Furthermore,
\[
\int_{t=0}^{\infty} \int_{u=0}^t \big(e^{-\beta (t-u)} M_{t-u} - n_1^{\bullet})\,e^{-\beta u}\,dG_y(u)\,dt =  V(y) \int_t^{\infty} \big(e^{-\beta t} M_t - n_1^{\bullet})\,dt\:,
\]
so that the integrability of $e^{-\beta t} M_t[y] - n_1^{\bullet} \mu V(y)$ follows from that of $e^{-\beta t} M_t - n_1^{\bullet}$. But the latter is just Proposition 8 in \cite{j}. \hfill $\Box$\\

Lemma \ref{statt} fails in the Galton-Watson scenario, and the reason for this is just the lack of integrability of $e^{-\beta t} M_t - n_1^{\bullet}$ in this case ($n_1^{\bullet}$ should be taken equal to the value of $\kappa$ in Equation \reff{triple}). This is not in contradiction with Proposition \ref{frrat}; it only shows that with the methods we use for the proof of, ultimately, Theorem \ref{main}, we cannot prove or disprove either of Propositions \ref{cthgw} and \ref{frrat}.

We now come to the main result of this section:
\begin{lemma} \label{thbh} Suppose that the life-time distribution $G^{\bullet}(t)$ of mutant cells is not lattice-like, and that $g(s)$ is as given in Equation \reff{g} with $\omega(s) = s$. Then, with $n^{\bullet}_1$ as defined in Equation \reff{n1}, we have that
\begin{displaymath}
\lim_{s \to 1^-} \frac{1 - g(s)}{1 - s} + n^{\bullet}_1 \log(1 - s)
\end{displaymath}
exists and is bounded away from $\pm\infty$.
\end{lemma}
{\it Proof}. First observe that $g(s)\equiv\syG[0](\fs)$, where $\fs(x)$ is {\it that} function on $[0, \infty)$ which assumes the constant value $s$ for every $x$, and $\syG[0](\fs)$ is as given in Lemma \ref{statt}. Therefore it will be sufficient to show that
\begin{displaymath}
\frac{1 - \syG[0]\big(R(e^{- \beta t})\big)}{\log R(e^{- \beta t})} + n^{\bullet}_1 \log\big(-\log R(e^{- \beta t})\big)
\end{displaymath}
converges as $t \to \infty$. By Lemma \ref{statt}, this boils down to proving convergence for
\begin{equation}
\frac{1 - \syG[0](\fs_t)}{\log R(e^{- \beta t})} + \log \big(-\log R(e^{- \beta t})\big) \:, \label{corr}
\end{equation}
and to demonstrate {\it this}, we define
\begin{equation}
\syF[x](\fr) :=: \syF_{\delta}[x](\fr) := \frac{1}{\delta'} \int_0^{\delta'} e^{- \beta u} \Phi_u[x](\fr)\,du \:
\end{equation}
and
\begin{equation}
\syH[x](\fr) :=: \syH_{\delta}[x](\fr)  \frac{\syF_{\delta}[x]\big(\fone\big) - \syF_{\delta}[x] \big(\Phi_{\delta'}[\cdot](\fr)\big)}{\syF_{\delta}[x](\fone) - \syF_{\delta}[x](\fr)}\:.
\end{equation}
We need to know about the behaviour of these two functionals in some neighbourhood of the unit function $\fone$. First, $\syF[x](\fone) = (1 - e^{-\delta})/\delta$. Now choose any $t > 0$ and observe that, because of Equation \reff{recur},
\begin{equation}
\frac{\syF[x](\fone) - \syF[x](\fs_t)}{\log R(e^{- \beta t})} = \frac{1}{\delta'} \int_0^{\delta'} e^{- \beta u} \, \frac{1 - \fs_{t-u}(x)}{\log R(e^{- \beta t})} \, du\to - \mu V(x) \label{nolog}
\end{equation}
as $t\to\infty$. Moreover, by Lemma \ref{straight}, the difference between the two sides is bounded in absolute value by some constant times $\mu - h\circ R(e^{- \beta (t-\delta')})$. Next,
\begin{eqnarray}
\lefteqn{\syH[x](\fs_t) = \frac{1}{\delta'} \int_0^{\delta'} e^{- \beta u} \, \frac{\Phi_u[x](\fone) - \Phi_{u + \delta'}[x](\fs_t)}{\syF[x](\fone) - \syF[x](\fs_t)} \, du} \label{hf} \\
&=& \frac{e^{\beta\delta'}}{\delta'} \int_0^{\delta'} e^{- \beta (u+\delta')} \, \frac{1 - \fs_{t - u - \delta'}(x)}{\log R(e^{- \beta t})}\,\frac{\log R(e^{- \beta t})}{\syF[x](\fone) - \syF[x](\fs_t)}\,du\to e^{\delta}\:,\nonumber
\end{eqnarray}
{\it iff mutation is neutral}, and convergence is at least as quick as $\sim \mu - h\circ R(e^{- \beta (t-2\delta')})$, again by Lemma \ref{straight}.

Returning now to the proof of Lemma \ref{thbh}, it is readily checked that, by the definition of the functional $\syG[x](\fr)$ in Lemma \ref{statt},
\begin{equation}
\syG[x](\fr) = \delta\sum_{i=0}^{\infty} \frac{\syF[x]\big(\Phi_{i\ast \delta'}[\cdot](\fr)\big)}{e^{i\delta}} \:.
\end{equation}
Thus, we obtain in complete analogy with Equation \reff{triplem},
\begin{eqnarray}
\lefteqn{\frac{1 - \syG[x](\fr)}{\log R(e^{-\beta t})} + \frac{\syF[x](\fone)-\syF[x](\fr)}{\log R(e^{-\beta t})}\,\log\big(\syF[x](\fone) - \syF[x](\fr)\big)} \label{triplebh} \\
&=&\delta\sum_{i=0}^{\infty} \frac{1}{e^{i\delta}} \, \frac{\syF[x](\fone)-\syF[x]\big(\Phi_{i\ast\delta'}[\cdot](\fr)\big)}{\log R(e^{-\beta t})}\nonumber \\
&-& \frac{1}{e^{\delta}}\sum_{i=0}^{\infty}\Bigg(\frac{\syH[x] \big(\Phi_{i\ast\delta'}[\cdot](\fr)\big)}{e^{i\delta}}\,\frac{\syF[x](\fone)-\syF[x]\big(\Phi_{i\ast\delta'}[\cdot](\fr)\big)}{\log R(e^{-\beta t})}\times\nonumber \\
&& \qquad \qquad \log \Big(\syH[x]\big(\Phi_{i\ast\delta'}[\cdot](\fr)\big) \Big) \Bigg) \nonumber \\
&+& \frac{1}{e^{\delta}}\sum_{i=0}^{\infty} \Bigg(\frac{e^{\delta} - \syH[x]\big(\Phi_{i\ast\delta'}[\cdot](\fr) \big)}{e^{i\delta}}\,\frac{\syF[x](\fone) - \syF[x]\big(\Phi_{i\ast\delta'}[\cdot](\fr)\big)}{\log R(e^{-\beta t})}\times \nonumber \\
&& \qquad \qquad \log\Big(\syF[x](\fone)-\syF[x]\big( \Phi_{i\ast\delta'}[\cdot](\fr)\big)\Big)\Bigg)\nonumber\:.
\end{eqnarray}
By a similar line of argument as for the proof of Lemma \ref{thctm}, we find that the expression on the right-hand side of the above equation, when taken at $\fr(\cdot) = \fs_t(\cdot)$, converges as $t\to\infty$ through integer multiples of $\delta'$. As an illustration, consider the final right-hand term in the above array, for which we claim that
\begin{eqnarray}
\lefteqn{\frac{1}{e^{\delta}}\sum_{i=0}^{\infty} \Bigg(\frac{e^{\delta} - \syH[x]\big(\Phi_{i\ast\delta'}[\cdot](\fs_0) \big)}{e^{i\delta}}\,\frac{\syF[x](\fone) - \syF[x]\big(\Phi_{i\ast\delta'}[\cdot](\fs_0)\big)}{e^{n \delta} \log R(e^{-n\delta})}}\nonumber \\
&& \qquad \qquad \times \log\Big(\syF[x](\fone)-\syF[x]\big( \Phi_{i\ast\delta'}[\cdot](\fs_0)\big)\Big)\Bigg)\nonumber\\
&&+\:\frac{1}{e^{\delta}}\sum_{i=1}^n e^{i\delta}\Bigg(\frac{e^{\delta} - \syH[x](\fs_{i\ast\delta'})}{e^{n\delta}}\, \frac{\syF[x](\fone) - \syF[x](\fs_{i\ast\delta'})}{\log R(e^{-n\delta})}\nonumber\\
&& \qquad \qquad \times \log\Big(\syF[x](\fone)-\syF[x](\fs_{i\ast\delta'})\Big)\Bigg)\nonumber
\end{eqnarray}
converges as $n\to\infty$. But this is clear from Equations \reff{nolog}, \reff{hf}, assumption ($\star$), and the fact that $e^{\beta t}\log R(e^{-\beta t})$ approaches a finite limit as $t$ tends to infinity \cite{sch}.

Thus we have shown that
\begin{eqnarray}
\lefteqn{\lim_{n\to\infty}\frac{1 - \syG[x](\fs_{n\delta'})}{\log R(e^{-n\delta})} + \frac{\syF_{\delta}[x](\fone)-\syF_{\delta}[x](\fs_{n\delta'})}{\log R(e^{-n\delta })}\log\big(-\log R(e^{-n\delta})\big)}\nonumber\\
&&+\:\frac{\syF_{\delta}[x](\fone)-\syF_{\delta}[x](\fs_{n\delta'})}{\log R(e^{-n\delta})}\log\left(\frac{\syF_{\delta}[x](\fone) - \syF_{\delta}[x](\fs_{n\delta'})}{-\log R(e^{-n\delta})}\right)\nonumber
\end{eqnarray}
exists for fixed $\delta$. We now have to deal with the fact that $\syF_{\delta}[\cdot](\cdot)$ depends on $\delta$ as an external parameter. Because of \reff{nolog}, the final term in the above expression has a finite limit which does not depend on $\delta$, and from the little discussion following \reff{nolog} it is clear that for $0<\tau<\delta$ the difference
\[
\frac{\syF_{\tau}[x](\fone)-\syF_{\tau}[x](\fs_t)}{\log R(e^{-\beta t)}} - \frac{\syF_{\delta}[x](\fone)-\syF_{\delta}[x](\fs_t)}{\log R(e^{-\beta t})}
\]
can be of order at worst $\mu - h\circ R(e^{-\beta(t-\delta)})$. It follows that
\[
\lim_{n\to\infty}\frac{1 - \syG[x](\fs_{n\delta'})}{\log R(e^{-n\delta})} + \frac{\syF_{\tau}[x](\fone)-\syF_{\tau}[x](\fs_{n\delta'})}{\log R(e^{-n\delta })}\log\big(-\log R(e^{-n\delta})\big)
\]
{\it also} exists for fixed $\tau$ and arbitrary $\delta > 0$
Now
\[
\syG[x](\fs_t) = \beta \int_0^t e^{-\beta u} \fs_{t-u}(x)\,du + e^{-\beta t} \beta \int_0^{\infty} e^{-\beta u} \Phi_u[x](\fs_0) \, du
\]
and
\[
\syF_{\tau}[x](\fs_t) = \frac{\beta e^{- \beta t}}{\tau} \int_{t-\tau/\beta}^t e^{\beta u} \fs_u\,du
\]
are continuous in $t$, and the same is obviously true for $R(e^{-\beta t})$. Thus it follows again by Kingman's lemma that
\[
\lim_{t\to\infty}\frac{1 - \syG[x](\fs_t)}{\log R(e^{-\beta t})} + \frac{\syF[x](\fone)-\syF[x](\fs_t)}{\log R(e^{-\beta t})}\log\big(-\log R(e^{-\beta t})\big)
\]
exists in $(-\infty,\infty)$, which by Equation \reff{nolog} implies that the same is true for
\[
\lim_{t\to\infty}\frac{1 - \syG[x](\fs_t)}{\log R(e^{-\beta t})} - \mu V(x) \log\big(-\log R(e^{-\beta t})\big)\:.
\]
But this proves the claim made for \reff{corr}, since $\mu V(0) = 1$, and thus the proof of Lemma \ref{thbh} and, finally, that of Theorem \ref{main} is complete.\hfill$\Box$
\section{Discussion}
The good news to take home from this paper is that, under the condition that $m$ is large enough for most cultures in a fluctuation experiment to produce a fairly large number of mutant colonies (by comparison with Lea and Coulson's \cite{lc} results, one would venture that $m \geq 4$ is already sufficient), it is possible to forget about the subtleties of clonal proliferation and still come up with an estimate for the mutation rate that makes better use of the data than just looking at the fraction of cultures without any mutants. Indeed, starting from Lea and Coulson's observation that for suitably chosen constants $a$, $b$, and $c$, the variate $a/(\xi + b) - c$ (with $\xi$ as in Theorem \ref{main}) is approximately normal distributed, one would go about evaluating the experiment as exemplified in \cite{lc}, or rather (because of the additional summand $\delta/n_1^{\bullet}$) as developed further in \cite{wa2}. The advantage of this method is that additional factors like phenotypic lag of mutants, or residual growth of non-mutant bacteria on selection medium, readily integrate into the calculations. The bad news is that the estimate for $m$ is consistently {\it wrong} by a factor $n_1^{\bullet}$. This confirms an observation of Oprea and Kepler \cite{ok}, who have found that their continuum approximation of the Luria-Delbr\"uck distribution can be adapted to realistic cell-cycle time distributions by a two-parameter generalisation. Although Oprea and Kepler are only concerned with the case $\mu = 2$, one can see that their factor $b/2$ corresponds to our $n_1^{\bullet}$, and it seems that the natural logarithm of their factor $c$ corresponds to our $\delta/n_1^{\bullet}$.

It is natural to ask about the values of $n_1^{\bullet}$ and $\delta$ for realistic models of cell proliferation. Since $\delta/n_1^{\bullet}$ can in principle be estimated from the experimental data \cite{wa2}, the question of practical relevance really is about the magnitude of $n_1^{\bullet}$. For Kendall's \cite{kend} multi-stage model of cell proliferation, for instance, with the bacteria lingering in each of $k$ stages for an exponentially distributed length of time, the life-time distributions are Gamma-like with densities
\begin{displaymath}
\frac{k^k \lambda^k}{\Gamma(k)} e^{- k \lambda t} t^{k - 1} \:,
\end{displaymath}
which implies $\beta = k \lambda (\mu^{k^{-1}} - 1)$, and
\begin{displaymath}
n_1^{\bullet} = \frac{\mu - 1}{k(\mu - \mu^{1 - k^{-1}})} \:,
\end{displaymath}
so that $n_1^{\bullet}$ approaches $1/k$ for {\it very} large values of $\mu$. For Rahn's \cite{rahn} life-time distribution density
\begin{displaymath}
\alpha k e^{- \alpha t} (1 - e^{- \alpha t})^{k - 1} \:,
\end{displaymath}
which obtains under the assumption that cells linger in each of in total $k$ stages for a time that is Gamma distributed with shape parameter 2, the respective parameters cannot be evaluated explicitly anymore (at least not unless $k \leq 4$), but it seems that $\beta \approx \alpha (k + 1) (\mu^{k^{-1}} - 1)/2$, and
\begin{displaymath}
n_1^{\bullet} = \frac{\alpha}{\beta} \, \frac{\mu - 1}{\mu} \left( \sum_{i = 1}^k \frac{1}{\beta / \alpha + i} \right)^{-1}\:.
\end{displaymath}
It is not difficult to prove via H\"older's inequality that for {\it every} life-time distribution of bacteria,
\begin{displaymath}
n_1^{\bullet} \geq \frac{\mu - 1}{\mu\log\mu}\:,
\end{displaymath}
but (by taking for $G^{\bullet}(\cdot)$ the Gamma distributions) it is also not difficult to see that $n_1^{\bullet}$ cannot be bounded from above. Still, it would seem that for realistic models of clonal proliferation, $n_1^{\bullet}$ varies in the range from 0.7 to 1. This would result in an underestimation of the mutation rate by, at worst, 70 percent, which isn't {\it very} much, but for experiments with a sufficiently large number of cultures (say, about 100) could differ significantly from an unbiased estimate.

The question about the magnitude of $\delta$ seems interesting in itself. The results of Oprea and Kepler suggest that at least for the cell-cycle time distributions investigated in their paper \cite{ok}, $\delta$ is positive and of about the same magnitude as $n_1^{\bullet}$. The proofs of Lemmas \ref{thctm} and \ref{thbh} partially confirm this. In fact, since in practice the offspring distribution of (mutant and non-mutant) bacteria has moments of {\it any} order, the $\omega$ in condition $(\star)$ can be taken as positive, so that the sums in Equations \reff{triplem} and \reff{triplebh} converge geometrically, and are about $1 - \mu^{-1}$ in magnitude. Finally, it seems possible to demonstrate some \lq stability' of the Luria-Delbr\"uck distribution with respect to different models of clonal proliferation also if the assumption of neutrality of the mutation is dropped, although, unsurprisingly, one ends up with stable limit laws of index other than 1. It is planned to present these results in a future publication.\\

{\bf Acknowledgment.} Part of the work presented in this paper comes from the author's PhD thesis (University of Vienna), for which he is happy to acknowledge financial support by the Austrian Science Fund (FWF), Project P14682-N05. Thanks are due to J\"urgen Steiner, Wolfgang L\"offelhardt, and all the people at the L\"offelhardt lab for their hospitality, and to Reinhard B\"urger for his continuous and unfailing support.

\section{Appendix}
Our goal in this appendix is to supply a proof of Lemma \ref{straight}. Because
\[
\frac{1}{\log s} + \frac{1}{1-s} = \frac{1}{2}\frac{1+\frac{2}{3}(1-s)+\ldots}{1 + \frac{1}{2}(1-s) + \ldots}\:,
\]
we may as well prove the corresponding statement for $1-R(e^{-\beta t})$ in the denominator instead of $-\log R(e^{-\beta t})$. We have
\begin{eqnarray}
\lefteqn{\mu V(x) - e^{-\beta u} \frac{1 - \fs_{t-u}(x)}{1 - R(e^{-\beta t})}}\nonumber\\
&=& \int_0^{\infty} \left(\mu e^{-\beta v} - e^{-\beta u}\,\frac{1 - f \circ R(e^{-\beta(t-u+ v)})}{1 - R(e^{-\beta t})}\right)\,dG_x(v)\nonumber\\
&=& \int_0^{\infty} \big(\mu - h\circ R(e^{-\beta(t-u+v)})\big)e^{-\beta v}\,dG_x(v)\nonumber\\
&&-\:\int_0^{\infty} h\circ R(e^{-\beta(t-u+v)})\frac{\cX(t-u+v) - \cX(t)}{\cX(t)}\,e^{-\beta v}\,dG_x(v)\:,\label{start}
\end{eqnarray}
where
\begin{equation}
\cX(t) := e^{\beta t}\big(1-R(e^{-\beta t})\big)\:,
\end{equation}
and $R(\cdot)$ is the Laplace transform from Equation \reff{r}. By monotonicity of $h(\cdot)$, the first term on the right-hand side of the above Equation \reff{start} is no larger than
\[
V(x)\big(\mu - h\circ R(e^{-\beta (t-u)})\big)\:.
\]
The problem is hence to obtain an estimate for the second term which, as $\cX(t)$ is increasing in $t$, consists itself of a positive and a negative part. It is however easy to see that neither of these can be larger than
\begin{equation}
\int_0^{\infty} h\circ R(e^{-\beta(t-u+v)})\frac{\cX(t-u+v) - \cX(t-u)}{\cX(t-u)}\,e^{-\beta v}\,dG_x(v)\:,\label{s2}
\end{equation}
so we may as well attempt to bound this. Furthermore, since it is only the ratio $\cX(t-u+v)\cX(t-u)^{-1}$ which enters Equation \reff{s2}, we may, for ease of writing, assume that $\cX(t-u) = 1$ for $0\leq t-u$ arbitrary, but fixed. We now deduce from Equation \reff{was} that
\begin{equation}
\cX(t) = \int_0^{\infty} h \circ R(e^{-\beta (t+u)})\cX(t + u)\,e^{-\beta u}\,dG(u) \:,\label{t}
\end{equation}
and then, if we write
\begin{equation}
G_{\beta}(t) := \mu \int_0^t e^{-\beta u}\,dG(u)\:,\label{aux}
\end{equation}
that
\[
\cX(t) = \int_t^{\infty} \frac{h \circ R(e^{-\beta u}) - \mu}{\mu} \cX(u)\, dG_{\beta}(u-t)+\int_t^{\infty} \cX(u)\, dG_{\beta}(u-t)\:.
\]
This is a relation which can be iterated: we use \reff{t} to express the second integral in the above equation as
\begin{eqnarray}
\lefteqn{\cX(t) - \int_t^{\infty} \frac{h \circ R(e^{-\beta u}) - \mu}{\mu} \cX(u)\, dG_{\beta}(u-t)} \nonumber \\
&=& \int_{u=t}^{\infty} \left(\int_{v=u}^{\infty} \frac{h \circ R(e^{-\beta v})}{\mu}\cX(v)\, dG_{\beta}(v - u) \right) \, dG_{\beta}(u - t) \nonumber \\
&=& \int_{v=t}^{\infty} \frac{h \circ R(e^{-\beta v})}{\mu}\cX(v) \, d\left(\int_{u=t}^v G_{\beta}(v - u) \, dG_{\beta}(u - t)\right) \nonumber \\
&=& \int_0^{\infty} \frac{h \circ R(e^{-\beta (t + v)})}{\mu}\cX(t + v)\, dG_{\beta}^{\ast 2}(v) \:, \nonumber
\end{eqnarray}
where $G_{\beta}^{\ast 2}(t)$ is for the two-fold convolution of $G_{\beta}(t)$ with itself, and proceed by induction:
\begin{eqnarray}
\lefteqn{\cX(t) = \int_0^{\infty} \cX(t+u)\,dG_{\beta}^{\ast n}(u)}\nonumber\\
&&-\:\int_0^{\infty} \frac{\mu - h \circ R(e^{-\beta (t + u)})}{\mu} \cX(t+u) \sum_{i = 1}^n \,dG_{\beta}^{\ast i}(u) \label{est}
\end{eqnarray}
for every natural $1 \leq n < \infty$. Now the \lq renewal function\rq\
\begin{equation}
U_{\beta}(t) := \sum_{i = 1}^{\infty} G_{\beta}^{\ast i}(t) \label{rn}
\end{equation}
(the standard definition would be to start the summation with $i = 0$) is finite on bounded intervals \cite{an}, whence it follows that
\[
\lim_{n\to\infty} \int_0^{\infty} \cX(t+u)\,dG_{\beta}^{\ast n}(u) = \lim_{u\to\infty}\cX(t+u) =: \cX_{\infty}\:,
\]
where the limit exists in $[0,\infty]$ because the other terms in Equation \reff{est} are obviously monotonous in $n$. Assumption ($\star$) implies that, in fact, $\cX_{\infty} < \infty$ \cite{sch}. We finally obtain
\begin{equation}
\cX_{\infty} - \cX(t) = \int_0^{\infty} \frac{\mu - h \circ R(e^{-\beta (t + u)})}{\mu}\,\cX(t+u)\,dU_{\beta}(u)\:.
\end{equation}
By renewal theory \cite{f},
\[
1 + U_{\beta}(t) = \nu_{\beta}t + \tilde{U}_{\beta}(t)
\]
for $t>0$, where
\[
\nu_{\beta}^{-1} := \int_0^{\infty} t\,dG_{\beta}(t) = \mu \int_0^{\infty} t e^{-\beta t}\,dG(t)\:,
\]
and $\tilde{U}(\cdot)$ is a measure which (is concentrated on the nonnegative reals and) integrates to
\[
\sigma_{\beta} := \frac{\nu_{\beta}^2\mu}{2} \int_0^{\infty} t^2 e^{-\beta t}\,dG(t)\:.
\]
Thus
\begin{eqnarray}
\cX_{\infty} - \cX(t) &=& \nu_{\beta}\int_t^{\infty} \frac{\mu - h \circ R(e^{-\beta u})}{\mu}\,\cX(u)\,du\nonumber\\
&&+\:\int_0^{\infty} \frac{\mu - h \circ R(e^{-\beta (t + u)})}{\mu}\,\cX(t+u)\,d\tilde{U}_{\beta}(u)\:,\nonumber
\end{eqnarray}
and because $h\circ R(e^{-\beta t}) < \mu$ for $t\in[0,\infty)$, it then follows that
\begin{eqnarray}
\lefteqn{\int_0^{\infty} h\circ R(e^{-\beta(t-u+v)})\big(\cX(t-u+v) - 1\big)\,e^{-\beta v}\,dG_x(v)}\nonumber\\
&<&\nu_{\beta}\int_{v=0}^{\infty} \int_{w=t-u}^{t-u+v} \big(\mu - h \circ R(e^{-\beta w})\big)\,\cX(w)\,dv\,e^{-\beta v}\,dG_x(v)\nonumber\\
&&+\:\int_{v=0}^{\infty} e^{-\beta v}\,dG_x(v)\int_{w=0}^{\infty} \big(\mu - h \circ R(e^{-\beta (t-u+w)})\big)\,\cX(t-u+w)\,d\tilde{U}_{\beta}(w)\nonumber\\
&<&\frac{\cX_{\infty}}{\cX(0)}\big(\mu - h \circ R(e^{-\beta (t-u)})\big)\left(\sigma_{\beta} V(x) + \nu_{\beta}\int_0^{\infty} u \,e^{-\beta u}\,dG_x(u)\right)\:,\nonumber
\end{eqnarray}
where we have divided by $\cX(0)$ to get rid of the assumption $\cX(t-u) = 1$. This obviously concludes the proof of the lemma, since $V(x)\leq 1$, and the integral is bounded by $(e\beta)^{-1}$.\hfill$\Box$

\end{document}